\documentclass[12pt, leqno]{article}
\usepackage{amsmath}
\usepackage{amssymb}
\usepackage{amsthm}

\def\underset#1#2{{\mathrel{\mathop {{}_{} {#2}}\limits_{{#1}_{}}}}}
\def\upplim_#1{\underset{#1}{\overline\lim}\;}
\def\lowlim_#1{\underset{#1}{\underline\lim}\;}
\newcommand{\C}{{\mathbf{C}}}
\newcommand{\codim}{{\mathrm{codim}}\,}

\newcommand{\G}{{\mathbf{G}}}
\newcommand{\I}{{\mathcal{I}}}
\renewcommand{\ker}{{\mathrm{Ker}\,}}
\newcommand{\N}{{\mathbf{N}}}

\newcommand{\R}{{\mathbf{R}}}

\newcommand{\st}{{\mathrm{Stab}}}
\newcommand{\supp}{{\mathrm{Supp}}}
\newcommand{\Z}{{\mathbf{Z}}}
\newcommand{\f}{{\mathbf{f}}}
\newcommand{\OS}{{\mathcal O}_S}
\newtheorem{cl}[equation]{{\it Claim}}
\newtheorem{cor}[equation]{Corollary}

\newtheorem{ex}[equation]{\bf Example}
\newtheorem{klem}[equation]{{\bf Key Lemma}}
\newtheorem{llist}[equation]{\rm  }
\newtheorem{thm}[equation]{\bf Theorem}
\newtheorem{prop}[equation]{\bf Proposition}
\newtheorem{lem}[equation]{\bf Lemma}
\numberwithin{equation}{section}
\newenvironment{pf}
{{\it Proof.  }\hskip10pt} {\hfill{\it Q.E.D.}\par\vskip+10pt}

\title{A New Unicity Theorem and Erd\"os' Problem for
Polarized Semi-Abelian Varieties
\thanks{Research supported in part by Grant-in-Aid
for Scientific Research (S) 17104001.}}
\date{}
\author{Pietro Corvaja and Junjiro Noguchi}
\begin{document}
\parindent12pt
\baselineskip18pt
\maketitle
\begin{abstract}
In 1988 P. Erd\"os asked if the prime divisors of $x^n -1$ for all $n=1,2, \ldots$
determine the given integer $x$; the problem was affirmatively
answered by  Corrales-Rodorig\'a\~nez and R. Schoof \cite{crs}
in 1997 together with its elliptic version. 
Analogously, K. Yamanoi \cite{y04} proved in 2004 that the support
of the pull-backed divisor $f^{*}D$ of an ample divisor
on an abelian variety $A$ by an algebraically non-degenerate
entire holomorphic curve $f: \C \to A$
essentially determines the pair $(A, D)$. 

By making use of the main theorem of 
\cite{nwy08} 
we here deal with this problem for semi-abelian varieties: namely, given two polarized semi-abelian varieties $(A_1, D_1)$, $(A_2,D_2)$ and entire non-degenerate holomorphic curves $f_i:\C\to A_i$, $i=1,2$, we classify the cases when  the inclusion $\supp f_1^*D_1\subset \supp f_2^* D_2$ holds.
We  also apply the main result of \cite{cz}
to prove an arithmetic counterpart.
\end{abstract}

\section{Introduction and main results.}

The purpose of this paper is to study a kind of unicity problem
for semi-abelian varieties with given divisors (polarization)
in terms of entire holomorphic curves and of arithmetic recurrences.

Let $A_i$ ($i=1,2$) be semi-abelian varieties and let
$D_i$  be reduced divisors on $A_i$.
Just for the sake of simplicity, we assume here that $D_i$ $(i=1,2)$
are irreducible and have trivial stabilizers  
$$
\st(D_i)=\{0\},\quad i=1,2,
$$
(see \S 2 for the notation).
Note that these assumptions are not restrictive by a reduction
argument: Cf.\ \S4 for the general case.

Our first result is as follows (cf.\ \S2 for the notation):
\smallskip

\begin{thm}
\label{main}
Let $f_i: \C \to A_i$ $(i=1,2)$ be non-degenerate holomorphic
curves.
\begin{enumerate}
\item
Assume that
\begin{align}
\label{1.1}
\underline{\supp f_1^*D_1}_\infty &\subset
\underline{\supp f_2^*D_2}_\infty , \\
\label{1.2}
N_1(r, f_1^*D_1) &\sim N_1(r, f_2^*D_2)||^{\empty}.
\end{align}
Then there is a (finite) \'etale morphism $\phi:A_1 \to A_2$ such that
$\phi \circ f_1=f_2$ and $D_1 \subset \phi^* D_2$.
\item
If equality holds in \eqref{1.1}, then
$\phi:A_1 \to A_2$ of {\rm (i)} is an isomorphism and $D_1 = \phi^* D_2$.
\item
If $A_i~ (i=1,2)$ are abelian varieties and $D_2$ is smooth or more
generally locally irreducible at every point of $D_2$,
then $\phi$ of {\rm (i)} is an isomorphism and $D_1 = \phi^* D_2$.
\end{enumerate}
\end{thm}

{\bf N.B.}  (i) When $A_i$ are abelian varieties $(i=1,2)$,
the above (ii) is K. Yamanoi's Unicity Theorem
(\cite{y04}, Theorem 3.2.1).

(ii) Because of the proof, assumption \eqref{1.1} can be replaced by
the following estimate:
$$
N_1(r, f_1^{*} D_1) - N_1(r, f_1^{*} D_1 \cap f_2^* D_2 )=
o(N_1(r, f_1^*D_1))||.
$$

(iii) Assumption \eqref{1.2} is necessary (cf. Example \ref{ce}).
\medskip

The following corollary follows immediately from Theorem \ref{main}.

\begin{cor}
{\rm (i)} Let $f_1: \C \to \C^*$ and $f_2:\C \to A_2$ with an elliptic curve
$A_2$ be holomorphic and non-constant.
Then $\underline{f_1^{-1}\{1\}}_\infty \not=
 \underline{f_2^{-1}\{0\}}_\infty$, where $0$ is the zero element
of $A_2$.

{\rm (ii)} If $\dim A_1 \not= \dim A_2$ in Theorem \ref{main}, then
$\underline{f_1^{-1}D_1}_\infty \not= \underline{f_2^{-1}D_2}_\infty$.
 
\end{cor}

{\bf N.B.} (i) The first statement means that the difference
of the value distribution
property caused by the quotient
$\C^{*} \to \C^{*}/\langle \tau \rangle=A_2$ can not be
recovered by the choice of $f_1: \C \to \C^{*}$ and $f_2: \C \to A_2$,
even though $f_i$ are
allowed to be arbitrarily transcendental.

(ii) The second statement implies that the distribution
of $f_i^{-1}D_i$ about
$\infty$ contains the topological information of $\dim A_i$;
it is interesting to observe
that this works even for one parameter subgroups with Zariski dense image.

\bigskip

Due to the well-known correspondence between Number Theory and
Nevanlinna Theory, it is tempting to give a number-theoretic analogue of
Theorem \ref{main}.
In 1988, P\'al Erd\"os raised the following problem:

Erd\"os' Problem.
{\it  Let $x, y$ be positive integers with the property that
for all positive integers $n$ the set of prime numbers dividing $x^n -1$
is equal to the set of prime numbers dividing $y^n -1$. Is then $x=y$? 
}

In \cite{crs} Corrales-Rodorig\'a\~nez and R. Schoof
solved this problem, by proving that if,  for every $n$, each prime
dividing $x^n-1$ divides also $y^n-1$, then $y$ is a power of $x$. They
also solved the natural elliptic analogue, for two (elliptic) recurrences in a same elliptic curve.

A related problem asks to classify the cases where $x^n-1$ divides
$y^n-1$; in \cite{czIndag},
it was proved that the same conclusion holds
(i.e. $y$ is a power of $x$) under the hypothesis
that $x^n-1$ divides $y^n-1$ for
infinitely many positive integers $n$. Note that the conclusion that $y$
is a power of $x$ can be translated in geometric terms by saying that
there exists an isogeny $\phi:\G_m\to\G_m$ such that
$\phi(x)=y$
and $\phi^* \{1\}\supset \{1\}$.  

The natural generalization to several variables is represented by Pisot's
problem, asking to characterize the pairs of linear recurrent sequences
$(n\mapsto \f_1(n)), (n\mapsto \f_2(n))$ such that $\f_1(n)$ divides
$\f_2(n)$ for every  integer  $n$ (or for infinitely many integers $n$).
We shall explain in \S 2 why  these problems can be viewed as
some analogue of the unicity problem for holomorphic maps to
semi-abelian varieties.
A first case (where it is assumed that the divisibility holds for
every integer $n$) was solved by van der Poorten in \cite{vdP} (see also \cite{rum}),
while  in \cite{cz} this was generalized to the case when the divisibility
is assumed to hold  just for infinitely many integers $n$.  

We would like to deal with the case of a semi-abelian variety
with a given divisor, i.e., a polarized semi-abelian variety.
As it often happens, the complex-analytic theory is more advanced, and
we dispose only of partial results in the arithmetic case. In the present
situation, we can prove an analogue of 
Theorem \ref{main} in the toric case,
but not in the general case of semi-abelian varieties,
that is left to be a {\it Conjecture}.

Here is our result in the arithmetic case;
 again, the notation is explained in \S2:

\begin{thm}\label{with_multiplicity}
Let $\OS$ be a ring of $S$-integers in a number field $k$.
Let $\G_1$ and $\G_2$ be linear tori, and let
$g_i\in\G_i(\OS)$  be elements generating Zariski-dense subgroups in
 $\G_i$ $(i=1,2)$.
Let $D_i$ be reduced divisors defined over $k$, with   defining ideals
$\I(D_i)$, such  that each  irreducible component
has a  finite  stabilizer and the stabilizer of $D_2$ is trivial.
 Suppose that for
infinitely many natural numbers $n\in\N$, the inclusion of ideals
\begin{equation}
\label{1.5}
 (g_1^n)^* \I(D_1)\supset (g_2^n)^* \I(D_2)
\end{equation}
holds.
Then there exists an \'etale morphism
$\phi: \G_1\to\G_2$, defined over $k$,
and a positive integer $h$ such that 
$\phi(g_1^h)=g_2^h$ and
$ D_1 \subset \phi^{*} (D_2)$.
\end{thm}
 
\noindent
(Note that inclusion of ideals \eqref{1.5} is a stronger version of the
inclusion \eqref{1.1}, although it apparently goes in the opposite
direction!) 

{\bf N.B.} (i) Theorem \ref{with_multiplicity} is essentially due to the
main results of \cite{cz}. We shall derive it from Corollary 1 of
\cite{cz}.

(ii) Example \ref{es} will show that we cannot take
$h=1$ in general. 

(iii) By Example \ref{es2}, the condition on the finiteness of the
stabilizer of each component of $D_1$ and $D_2$  cannot be omitted, nor can
the condition on the triviality of the stabilizer of $D_2$. However, it
is easy to derive from our proof a nontrivial consequence even in the
general case, see Theorem \ref{gene}.

(iv) We assume that inequality (inclusion) \eqref{1.5} of ideals
holds only for an infinite sequence of $n$,
 not necessarily for all large $n$.
On the contrary, we need the inequality of ideals,
not only of their {\it  supports}, i.e. of
the primes containing the corresponding ideals.
\smallskip

One might ask for a similar conclusion assuming only the inequality of
supports,
\begin{equation}
\label{1.7}
\supp\, (g_1^n)^* \I (D_1)\subset  \supp\, (g_2^n)^* \I(D_2),
\end{equation}
but we then need \eqref{1.7} to hold for {\it every}  exponent
$n$ (this is in analogy with Erd\"os' original problem).
In this case, using a work of Barski, B\'ezivin and Schinzel \cite{bbs},
we will show some result which is slightly weaker
than the analogue of the above theorems
(see Proposition \ref{no_multiplicity}).

{\it Acknowledgement.}  The present joint research was initiated
when the authors shared a time at the Fields Institute, Toronto,
participating in the Thematic Program on Arithmetic Geometry,
 Hyperbolic Geometry and Related Topics, July-December, 2008.
The second author had a useful discussion on the difference of
divisors on abelian and semi-abelian varieties with our colleague,
Professor J\"org Winkelmann.
We would like to express our sincere gratitude to the
Fields Institute, and Professor J. Winkelmann.

\section{Notation.}

In this paragraph we introduce the necessary notation, both in the
analytic setting (a), and in the arithmetic one (b).
\medskip

(a) A subset $Z \subset \C$ determines a germ, denoted by
$\underline{Z}_\infty$, of subsets at the infinity of
the Riemann sphere $\C \cup\{\infty\}$.
For two subsets $Z_i \subset \C, i=1,2$, the relation
$$
\underline{Z_{1}}_\infty \supset \underline{Z_2}_\infty
$$
makes sense: it means there is a number $r_0>0$ such that
$$
Z_1 \cap \{z\in \C; |z|>r_0\} \supset Z_2 \cap \{z\in \C; |z|>r_0\}.
$$

For two functions $\phi_i(r)\geq 0,
r >0, i=1,2$, we write
$$
\phi_1(r) \leq \phi_2(r)||,
$$
if there is an exceptional (Borel) subset $E\subset (0, \infty)$
with a finite measure $m(E) < \infty$, satisfying
$$
\phi_1(r) \leq \phi_2(r),\quad r \in (0,\infty)\setminus E.
$$

If there is a positive constant $C>0$ satisfying
$\phi_1(r) \leq C\phi_2(r)||$, we write
$$
\phi_1(r)=O(\phi_2(r))||.
$$
If $\phi_1(r)\leq \epsilon \phi_2(r)||$ for an arbitrary $\epsilon >0$,
we write
$$
\phi_1(r)=o(\phi_2(r))||.
$$
Here, it is noted that the exceptional subset in
$\phi_1(r)\leq \epsilon \phi_2(r)||$ may depend on $\epsilon>0$.
We write
\begin{equation}
\label{equiv}
\phi_1(r) \sim \phi_2(r)||,
\end{equation}
if $\phi_1(r)=O(\phi_2(r))||$ and $\phi_2(r)=O(\phi_1(r))||$.
If the exceptional subset of $(0, \infty)$ is empty in the
above expressions, we will simply drop the symbol ``$||$''.

A complex algebraic group $A$ admitting a representation
$$
0 \to \G_m^t \to A \to A_0 \to 0,
$$
where $\G_m=\C\setminus\{0\}$ is the multiplicative group
and $A_0$ is an abelian variety,
is called a {\it semi-abelian variety}.
Let $A$ be a semi-abelian variety and let $B$ be an algebraic subset
of $A$.
Set
$$
\st(B)=\{x \in A; x+B=B\},
$$
which is called the {\it stabilizer} of $B$ (in $A$).
We denote by $\st(B)^0$ the identity component of $\st(B)$.
For a given $B$ there is an equivariant compactification
$\bar{A}$ (smooth) of $A$ such that the closure $\bar{B}$ of $B$ in
$\bar{A}$ contains no $A$-orbit (cf.\ \cite{nwy08}).
If $B$ is a reduced divisor $D$ on $A$, then
\begin{llist}
\qquad
\label{big}
such $\bar D$ is big on $\bar A$
if $\st(D)$ is finite {\rm (cf.\ \cite{nwy08})}.
\end{llist}

Let $f: \C \to A$ be a holomorphic curve.
If the image $f(\C)$ is (resp.\ is not) Zariski dense in $A$,
$f$ is said to be algebraically non-degenerate (resp.\ degenerate);
from now on,  we simply say that
$f$ is {\it non-degenerate} (resp.\ degenerate).
We use $T(r, \omega_{B, f})$ and $T_f(r, c_1(\bar{D}))$ for
the order functions,
$N(r, f^*D)$ and $N_k(r, f^*D)$ for the counting functions,
as defined as in \cite{no84}, \cite{nwy02} and \cite{nwy08}.

Let $f:\C \to M$ be a non-constant holomorphic curve into a
projective algebraic manifold $M$,
and let $T_f(r)$ denote the order function of $f$ with respect
to an ample divisor on $M$.
Note that there are several ways to define the order function of
$f$, cf.\ \cite{no84}; they are equivalent in the sense of \eqref{equiv}.
Then $f$ is rational if and only if $\overline{\lim}~ (\log r)/T_f(r) >0$
(cf.\ \cite{no84}).
Now, let $M$ be a compactification of a semi-abelian variety $A$,
and suppose that $f(\C) \subset A$.
Since there is no non-constant rational map from $\C$ into $A$, we have
$$
\log r =o(T_f(r)).
$$
Let $Z_i~(i=1,2)$ be effective divisors on $\C$ such that
$\underline{Z_1}_\infty = \underline{Z_2}_\infty$.
Then,
\begin{equation}
\label{logdiff}
N_k (r, Z_1)=N_k (r, Z_2 )+O(\log r)=
N_k (r, Z_2 ) +o(T_f(r))\quad (1 \leq k \leq \infty).
\end{equation}

(b) We now switch to the arithmetic setting. Let $k$ be a number field,
let $S$ be a finite set of places of $k$,
containing the archimedean ones,
and let $\OS$ be the corresponding ring of $S$-integers.
Let $\G\cong\G_m^n$ be a (split) linear torus; since the algebraic
variety $\G_m$ is canonically defined over the ring of rational
integers, we can tacitly considered it as a scheme $\G_{\OS}\to{\rm
Spec}(\OS)$.
Every integral point $g\in\G(\OS)$ can be viewed as a
morphism $g:{\rm Spec}(\OS)\to\G_{\OS}$. A divisor $D$ in $\G$,
defined over $\OS$, corresponds to an ideal $\I(D)$ of $\OS[\G]$;
its  pull-back $g^* \I(D)$ is naturally defined,
and is an ideal of $\OS$.
In fact, if $F(X_1,\ldots,X_n)=0$ is an equation for $D$,
where $F(X_1,\ldots,X_n)\in\OS[X_1,\ldots,X_n]$ is a polynomial,
then $g^*\I( D)$ is the ideal of $\OS$ generated by $F(g)$.
The support of this ideal, denoted by $\supp\, g^*(D)$,
is the set of maximal ideals containing $g^* \I(D)$.

\section{Proof of Theorem \ref{main}.}

Before beginning with the proof, we give a key lemma.
In general, let $f_i:\C \to A_i$ be non-degenerate holomorphic curves
into semi-abelian varieties $A_i$ ($i=1,2$).
Put
$
g=(f_1, f_2): \C \to A_1 \times A_2.
$
Let $A_0$ be the Zariski closure of $g(\C)$.
Then, by  Log Bloch-Ochiai's Theorem \cite{n81} $A_0$ is a translate
of a semi-abelian subvariety of $A_1 \times A_2$, so that
we have a non-degenerate holomorphic curve and the natural projections,
\begin{align}
\label{g}
&g: \C \to A_0,\\
\nonumber
& p_i : A_0 \to A_i\quad (i=1,2).
\end{align}
Since $f_i$ are non-degenerate, $p_i(A_0)=A_i$  ($i=1,2$).

\begin{klem}
\label{keylem}
Let the notation be as above.
Let $D_i$ be reduced divisors on $A_i$ for $i=1,2$.
Assume that
\begin{enumerate}
\item
$D_1$ is irreducible,
\item
$\underline{\supp f_1^*D_1}_\infty \subset
\underline{\supp f_2^*D_2}_\infty $,
\item
$N_1(r, f_1^*D_1) \sim N_1(r, f_2^*D_2)||$.
\end{enumerate}
Then, $p_1^*D_1 \subset p_2^*D_2$ on $A_0$.
\end{klem}

\begin{pf}
We reduce the case to the one where
\begin{equation}
\label{red}
\st(D_i)=\{0\}, \qquad i=1,2.
\end{equation}
For the reduction we set
\begin{align*}
&q_i : A_i \to A_i/\st(D_i)=A_i',\\
&q_i(D_i)=D_i/\st(D_i)=D_i',\\
&f_i'=q_i \circ f_i,\\
&g'=(f_1', f_2'): \C \to A_0'\: (\subset A_1' \times A_2')
\quad ( A_0' \hbox{ is the Zariski closure of } g'(\C)), \\
&p_i':A_0' \to A_i'\quad (\hbox{the natural projections}),\\
&\tilde{q}=(q_1,q_2):A_1 \times A_2 \to A_1' \times A_2',\\
&\tilde{q}_0=\tilde{q}|_{A_0}: A_0 \to A_0'.
\end{align*}
Then we see that $\st(D_i')=\{0\}$ and assumptions (i)$\sim$(iii) are
satisfied for $f_i'$ and $D_i'$.
Suppose that the present lemma was proved in this case.
Then we would have that ${p_1'}^* D_1' \subset {p_2'}^*D_2'$.
It follows that $\tilde{q}_0^* ( {p_1'}^* D_1') \subset
 \tilde{q}_0^* ({p_2'}^*D_2')$,
and hence that $ p_1^*(q_1^* D_1') \subset p_2^* (q_2^* D_2')$.
Therefore, $p_1^*D_1 \subset p_2^*D_2$.

Now we assume \eqref{red}.
By virtue of the second main theorem established by
\cite{nwy08} there exists for each $i=1,2$ an equivariant
compactification $\bar{A}_i$ of $A_i$
such that
\begin{equation*}
N_1(r, f_i^* D_i)=(1+o(1))T_{f_i}(r, c_1(\bar{D}_i))||,
\end{equation*}
where $\bar{D}_i$ is the closure of $D_i$.
Since $\bar{D}_i$ is big on $\bar{A}_i$ (cf.\ \ref{big}),
we may take the order function of $f_i$ by
$$
T_{f_i}(r)=T_{f_i}(r, c_1(\bar{D}_i)) \quad (i=1,2),
$$
and for $g=(f_1,f_2)$ by
$$
T_g(r)=T_{f_1}(r)+T_{f_2}(r).
$$
It follows that
$$
T_g(r) \sim T_{f_i}(r)|| \quad (i=1,2).
$$

Let $F$ be an arbitrary irreducible component of $p_1^*D_1$.
We are going to show

\begin{cl}
\label{claim1}
\quad
$F \subset p_2^*D_2$.
\end{cl}
There are finitely many elements $a_\nu \in \ker p_1, 1 \leq \nu \leq \nu_0$,
such that
$$
p_1^*D_1=\sum_{\nu=1}^{\nu_0} (F+a_\nu).
$$
Let $\bar{A}_0$ be an equivariant compactification of $A_0$ such that
the closure $\overline{p_1^*D}_1$ of $p_1^*D_1$ in $\bar{A}_0$
contains no $A_0$-orbit,
and $p_1$ extends to a morphism $\bar{p}_1: \bar{A}_0 \to \bar{A}_1$.
Then $c_1(\bar{F}+a_\nu)=c_1(\bar{F})$, and
\begin{align*}
T_{f_1}(r, c_1(\bar{D}_1)) &= T_g(r, c_1(\overline{p_1^*D}_1))
= \sum_{\nu=1}^{\nu_0} T_g(r, c_1(\bar{F}+a_\nu))\\
&=\nu_0 T_g(r, c_1(\bar{F})).
\end{align*}
Again, by the Second Main Theorem \cite{nwy08} and the above we have
$$
N_1(r, g^*F)=(1+o(1))T_g(r, c_1(\bar{F}))||\sim T_{f_1}(r)||.
$$
It follows from assumption (ii) and \eqref{logdiff} that
$$
N_1(r, g^*F)=N_1(r, g^* (F \cap p_2^*D_2))+O(\log r).
$$

Suppose that Claim \ref{claim1} is not true.
Then we would have
$$
\codim (F \cap \supp\, p_2^* D_2) \geq 2.
$$
The Main Theorem (ii) of \cite{nwy08} implies that
$$
N_1(r, g^* (F \cap p_2^*D_2))=o(T_g(r))||.
$$
Thus we obtain a contradictory estimate,
$$
T_{f_1}(r)=o(T_{f_1}(r))||.
$$
\end{pf}

{\it  Proof of Theorem \ref{main}.}
Firstly by \cite{n98} there are points $\zeta_i \in \C$ $(i=1,2)$
such that $f_i(\zeta_i) \in
D_i$. Then, considering the composites with the translations
$$
z \in \C \to f_i(z+\zeta_i)-f(\zeta_i),\quad D_i - f(\zeta_i),
$$
we may assume that $f_i(0)=0 \in D_i$ ($i=1,2$).
In what follows, we keep this setup.

(i) Let $g=(f_1,f_2): \C \to A_0 (\subset A_1 \times A_2)$ and 
$p_i:A_0 \to A_i$ be as above.
Note that $A_0$ is a semi-abelian subvariety of $A_1 \times A_2$.
Set $E_i=p_i^* D_i$. Let $F_{i}$ be an irreducible component of $E_i$
containing $0$.
Since $D_i$ are assumed to be irreducible,
there are finitely many points
$a_{i \nu} \in \ker p_i, 1 \leq \nu \leq \nu_i$
($a_{i1}=0$) such that
$$
E_i=\sum_{\nu=1}^{\nu_i} (a_{i\nu} + F_{i}) \quad (i=1,2),
$$
and $a_{i\nu} + F_{i}, 1 \leq \nu \leq \nu_i$, are mutually distinct for
each $i$.
By the Key Lemma \ref{keylem} we have in fact that
$$
F_1=F_2, \qquad E_1 \subset E_2.
$$
Note that $\st(E_1)^0=\st(E_2)^0=\st(F_1)^0\: (=\st(F_2)^0)$.
If $\dim \st(F_1)>0$, there should be a non-zero holomorphic vector field
$v$ on $A_0$ that is tangent to $E_1$ and $E_2$.
Therefore, the push-forward $p_{i*}v$ are tangent to $D_i$ ($i=1,2$).
Since $\st(D_i)=\{0\}$, $p_{i*}v=0$ ($i=1,2$), and hence $v=0$;
this is absurd.
Therefore we see that $\ker p_i$ are finite, and that
$$
\ker p_i=\{a_{i\nu}\}_{\nu=1}^{\nu_i}=\st(E_i).
$$
Since $F_1+a_{1 \nu}$ is an irreducible component of $p_2^*D_2$
and $D_2$ is irreducible, $p_2(F_1+a_{1 \nu})=
p_2(F_1)+p_2(a_{1\nu})=D_2+ p_2(a_{{1\nu}})=D_2$.
Since $\st(D_2)=\{0\}$, $p_2(a_{1\nu})=0$, so that
$a_{1 \nu} \in \ker p_2$. Therefore, $\ker p_1 \subset \ker p_2$,
and we have an isogeny $\phi:A_1 \to A_2$ by the composition of
the sequence of morphisms,
\begin{equation}
\label{isog}
A_1 \cong A_0/\ker p_1 \to A_0/\ker p_2 \cong A_2.
\end{equation}
It is immediate that $D_1 \subset \phi^*D_2$.

(ii)  Assume that equality holds in \eqref{1.1}.
The above proof implies that $E_1=E_2$ and $\ker p_1=\ker p_2$.
It follows from \eqref{isog}  that $\phi$ is an isomorphism.

(iii) It follows from the proof of (i) that $p_i : A_0 \to A_i$
($i=1,2$) are isogenies.
Since $\st(E_i)$ are finite, every irreducible component of $E_i$
is ample on $A_0$.
If there were two irreducible components $F'$ and $F''$ in $E_2$, 
$F' \cap F'' \not= \emptyset$.
Since $p_2: A_0 \to A_2$ is \'etale and $p_2(F')=p_2(F'')=D_2$,
$D_2$ is not locally irreducible; this is a contradiction.
Thus, $E_2$ is irreducible; since $E_1 \subset E_2$, $E_1=E_2$.
Hence, $\phi:A_1 \to A_2$ is an isomorphism.  \hfill {\it Q.E.D.}

\begin{ex}
\label{ce}
{\rm
Set $A_1=\C/\Z \:(\cong \G_m)$ and let $D_1=0$ be the zero element of
$A_1$.
Let $f_1:\C \to A_1$ be the covering map.
We take a number $\tau \in \C$ with the imaginary part $\Im \tau \not=0$.
Set $A_2=\C/(\Z+\Z \tau)$ which is an elliptic curve, and let
$D_2=0$ be the zero element of $A_2$.
Let $f_2: \C \to A_2$ be the covering map.

Then $f_1^{-1}D_1=\Z  \subset \Z+ \tau\Z= f_2^{-1} D_2$:
 assumption \eqref{1.1} of 
Theorem \ref{main} is satisfied.
There is, however, no non-constant
morphism $\phi:A_1 \to A_2$.
Note that
$$
N_1(r, f_1^*D_1)\sim r, \quad N_1(r, f_2^*D_2)\sim r^2.
$$
Thus, $N_1(r, f_1^*D_1) \not\sim N_1(r, f_2^*D_2)||$:
assumption \eqref{1.2} is failing.
}
\end{ex}

\section{The case of general $D_i$.}

Here we deal with the case where $D_i$ may be reducible
in 
Theorem \ref{main}.
The following example suggests that there must be some restrictions
on the given divisors $D_i$.

\begin{ex} {\rm
\label{ctex}
Let $\Z[i]$ denote the lattice of Gauss integers, and set
$$
\begin{array}{ll}
A_1=\C/\Z[i] \times \C^*, &\quad D_1
=\{0\}\times \C^*+\C/\Z[i]\times\{1\},\\
A_2=\C/\Z[i] \times \C/\Z[i], &\quad D_2=\{0\}\times \C/\Z[i]
+\C/\Z[i]\times\{0\}.
\end{array}
$$
Taking an irrational number $\alpha \in \R$, we set
\begin{align*}
&f_1: z \in \C \to ([z], e^{2\pi\alpha z}) \in A_1,\\
&f_2: z \in \C \to ([z], [\alpha z]) \in A_2,
\end{align*}
where $[z]$ denotes the point of $\C/\Z[i]$ represented by $z$.
Then $f_i$ are non-degenerate, $f_1^*D_1 \subset f_2^*D_2$,
and by calculation
$N_1(r, f_i^*D_i)\sim r^2$ ($i=1,2$).
There is, however, no morphism $\phi:A_1 \to A_2$ such that
$D_1 \subset \phi^*D_2$.
}
\end{ex}

The above example suggests that the stabilizers of the irreducible
components of $D_i$ should be restricted, while the assumption
of the triviality of $\st(D_i)$ is just a matter of reduction
to state the result as seen in \eqref{red}.

Let $A_i$ ($i=1,2$) be a semi-abelian variety and let $D_i$ be
reduced divisor on $A_i$, which may be reducible.
In the sequel we again assume that
\begin{equation}
\label{st}
\st(D_i)=\{0\}, \qquad i=1,2.
\end{equation}

\begin{thm}
\label{general}
Let $f_i: \C \to A_i$ ($i=1,2$) be a non-degenerate holomorphic curve.

{\rm (a)} Assume that
\begin{enumerate}
\item
every irreducible component of $D_1$ has a finite stabilizer;
\item
$\underline{\supp f_1^*D_1}_\infty \subset
\underline{\supp f_2^*D_2}_\infty$,
\item
$N_1(r, f_1^*D_1) \sim N_1(r, f_2^*D_2)||$.
\end{enumerate}
Then there exist a semi-abelian variety $A_0$, a non-degenerate
 holomorphic curve $g:\C \to A_0$,
reduced divisors $E_i$ on $A_i$ ($i=1,2$), and morphisms
$\phi_i:A_0 \to A_i$ such that
\begin{itemize}
\item
$E_1 \subset E_2$,
\item
$\phi_i^*D_i=E_i$ ($i=1,2$),
\item
$f_i=\phi_i \circ g \quad (i=1,2)$,
\item
$\phi_2: A_0 \to A_2$ is an \'etale morphism.
\end{itemize}

{\rm (b)} Moreover, if every irreducible component of $D_2$ has
a finite stabilizer, and
 equality holds in the above {\rm (ii)}, then
$\phi_i: A_0 \to A_i$ are isomorphisms, and $E_1=E_2$.
\end{thm}

\begin{pf}
(a) Let $g:\C \to A_0$ and $p_i: A_0 \to A_i$ be as in \eqref{g}.
Set $E_i=p_i^*D_i$.
We use the order functions $T_{f_i}(r)$ and $T_g(r)$ in the same sense
as in the proof of the Key Lemma \ref{keylem}.

Let $F$ be an arbitrary irreducible component of $E_1$.
Then $p_1(F)$ is an irreducible component of $D_1$, and by
assumption (i) $\st(p_1(F))$ is finite.
Hence, the closure $\overline{p_1(F)}$ in a compactification $\bar{A}_1$
of $A_1$ is big.
The second main theorem \cite{nwy08} implies that
$$
N_1(r, f_1^* (p_1(F)))\sim T_{f_1}(r)\sim N_1(r, f_1^*D_1)
\sim N_1(r, f_2^*D_2)||.
$$
We infer from the Key Lemma \ref{keylem} that
$F \subset E_2$; henceforth, $E_1 \subset E_2$.

By the same vector field argument as in the proof of
Theorem \ref{main} (i) we see that $\dim \st(E_2)=0$.
Thus $p_2:A_0 \to A_2$ is an \'etale morphism.
Setting $\phi_i=p_i$, we finish the proof of (a).

(b) In the same way as in (a) it immediately
follows from the assumptions that
$E_1=E_2$. Since $A_i \cong A_0/\st(E_i)$ ($i=1,2$), $\phi_i$
are isomorphisms.
\end{pf}

We have the following by Theorem \ref{general}.

\begin{cor}
Let $f_i:\C \to A_i$ be non-degenerate holomorphic curves
and let $D_i$ be reduced divisors on $A_i$ such that all irreducible
components of $D_i$ have finite stabilizers $(i=1,2)$.
If $A_1$ and $A_2$ are not isogenous,
then $\underline{f_1^{-1} D_1}_\infty \not=
 \underline{f_2^{-1} D_2}_\infty$. 
\end{cor}

\begin{pf}
By the assumption, $\st(D_i)$ are finite ($i=1,2$).
Then we have isogenies $q_i:A_i \to A_i/\st(D_i)$.
By setting $D_i'=q_i(D_i)$, the case is reduced to the one where
\eqref{st} is satisfied.
If $\underline{f_1^{-1} D_1}_\infty =
 \underline{f_2^{-1} D_2}_\infty$,
one would infer from Theorem \ref{general} that $\phi_i:A_0 \to A_i$ are
both \'etale morphisms for $i=1,2$, so that $A_1$ and $A_2$ are isogenous.
\end{pf}

\section{Proof of Theorem \ref{with_multiplicity}}

Theorem \ref{with_multiplicity} will be reduced to
statements about diophantine equations
involving linear recurrence sequences. For an introduction to the
general theory see \cite{vdP} or \cite{S}.
We shall actually be interested
in recurrence sequences of the form 
\begin{equation}\label{powersum}
 n\mapsto \f(n)=\sum_{i=1}^k
b_i\alpha_i^n, 
\end{equation}
where $k\geq 1$ is a natural number, $b_1,\ldots,b_k$
are nonzero complex numbers 
and $\alpha_1,\ldots,\alpha_k$ are nonzero pairwise distinct complex
numbers. The   representation \eqref{powersum} is unique: it  suffices
to note that the right-hand side in \eqref{powersum} cannot vanish for
$k$ consecutive values of $n$ (e.g. for $n=0,\ldots,k-1$) since the van
der Mond matrix $(\alpha_i^{n-1})_{1\leq i,n\leq k-1}$ is non-singular.
These recurrence sequences as in \eqref{powersum} will
be called {\it power sums}.
The complex numbers $\alpha_1,\ldots
\alpha_k$ will be called the {\it roots} of the power sum $\f$.
Theorem \ref{with_multiplicity} will be reduced to the following result,
appearing in a slightly more general formulation as Corollary 2 in
\cite{cz}:

\begin{lem}
\label{cz}
Let, for $i=1,2$,  $n\mapsto {\mathbf f}_i(n)$ be two power sums with
values in a ring of $S$-integers $\OS$. Suppose that the roots of
$\f_1,\f_2$ together generate a torsion-free multiplicative group. If
the ratio $\f_2(n)/\f_1(n)$ lies in $\OS$ for infinitely many $n$, then
the function $n\mapsto \f_2(n)/\f_1(n)$ is a power sum.
\end{lem}

Before starting the proof of Theorem \ref{with_multiplicity} we recall
some basic facts about the algebraic theory of power sums.
\medskip

Let $\mathcal{U}\subset\C^*$ be a torsion-free finitely generated
multiplicative group.
As an abstract group, $\mathcal{U}$ is isomorphic to $\Z^r$, where $r$ is
the rank of $\mathcal{U}$.
We shall be interested in the algebra of power
sums with roots in $\mathcal{U}$. Such an algebra is isomorphic to the
algebra $\C[\mathcal{U}]$, which in turn is isomorphic to the algebra of
Laurent polynomials
$\C[X_1,\ldots,X_r,X_1^{-1},\ldots,X_r^{-1}]=\C[\G_m^r]$.  Letting
$u_1,\ldots,u_r$ be generators of $\mathcal{U}$, an isomorphism is obtained
by sending the function $(n\mapsto u_i^n)$ to the monomial $X_i$ and
extending to a ring homomorphism. Note that the units $\C[X_1^{\pm
1},\ldots,X_r^{\pm 1}]$ are the monomials, so the units in the ring of
power sums are of the form $(n\mapsto a\alpha^n)$, for nonzero
$a\in\C^*$ and $\alpha\in\mathcal{U}$.

Let now $g\in\G_m^r$ be an element in a torus, not contained in any
subtorus. This means that if we write $g=(u_1,\ldots,u_r)$, the non-zero
complex numbers $u_1,\ldots,u_r$  are multiplicative independent,
i.e. generate a subgroup of maximal rank, which is necessarily
torsion-free.
The main link between the theory of linear recurrences and linear tori
is represented by the following fact: {\it For every regular function
$F\in\C[\G_m^r]$, the function $(n\mapsto F(g^n))$ is  a power sum with
roots in the group
$\mathcal{U}:=\langle u_1,\ldots,u_r \rangle$}. 

Another fact will be  repeatedly used in the sequel:

{\it Let
$F\in\C[\G_m^r]$ be a non zero regular function, $D=F^{-1}(0)$ the
corresponding divisor in $\G_m^r$. Then the stabilizer of $D$ has
positive dimension if and only if $F$ can be written, after applying an
automorphism to $\G_m^r$, in the form $X_r^l G(X_1,\ldots,X_{r-1})$, for
some integer $l$ (possibly zero) and a polynomial
$G(X_1,\ldots,X_{r-1})$ in $r-1$ indeterminates}.

Another equivalent
formulation is that the recurrence $(n\mapsto F(g^n))$ is of the form
$\alpha^n\f(n)$ where the roots of $\f$ generate a multiplicative group
of rank $<r$.

Let us now begin the proof of Theorem \ref{with_multiplicity}. Let
$d_i$, for $i=1,2$, be the dimension of the torus $G_i$. Let
$F_i(X_1,\ldots,X_{d_i})\in\C[X_1,\ldots,X_{d_i}]$ be Laurent
polynomials defining the irreducible divisor $D_i$:
$D_i=F_i^{-1}(0)$. We can clearly suppose that they are polynomials in
$X_1,\ldots,X_r$, and also that they have no monomial factor of positive
degree: both facts follow from the remark that by multiplying $F_i$ by a
monomial, the zero locus in $\G_m^{d_i}$ does not change.

Let $g_1,g_2$ be the elements  appearing in the statement of Theorem
\ref{with_multiplicity} and let $k\geq 1$ be the order of the torsion
subgroup of the group generated by the coordinates of
$g_1,g_2$.
Considering the partition of the set of natural numbers in
classes modulo $k$,   we obtain that for at least one integer
$r\in\{0,\ldots,k-1\}$   there will exist infinitely many positive
integers $m$ such that
$(g_1^{r+km})^* \I (D_1)\supset (g_2^{r+km})^* \I (D_2)$.
Replacing $g_i$ by $g_i^k$ and $D_i$ by its image
under the map $x\mapsto g_i^{-r}x$, we reduce the case to $k=1$. Note
that the conclusion we want to prove is not affected by this
replacement. So, from now on, we shall suppose that the   coordinates of
$g_1,g_2$ together generate a torsion-free multiplicative group
$\mathcal{U} = \langle u_1,\ldots,u_r \rangle \subset\C^*$.
Finally, let $\f_i(n)=F_i(g_i^n)$, so that $\f_1,\f_2$ are power sums
with roots in $\mathcal{U}$. Their values belong to the ring $\OS$,
although the roots and coefficients expressing $\f_1,\f_2$ as power sums
are not necessarily in the ring $\OS$.

The hypothesis of Theorem \ref{with_multiplicity} can be reformulated by
saying that for infinitely many natural numbers $n$, the ideal generated
in $\OS$ by $\f_1(n)$ contains the ideal generated by $\f_2(n)$,
i.e. $\f_1(n)$ divides $\f_2(n)$, in the ring $\OS$. Applying Lemma
\ref{cz}, we obtain that the power sum $\f_1$ divides the power sum
$\f_2$ in the ring of power sums. We shall exploit this fact, and see
how it leads to the sought conclusion of Theorem
\ref{with_multiplicity}.

\begin{lem}
Let $i\in\{1,2\}$. Suppose that the irreducible divisor
$D_i=F_i^{-1}(0)$ has finite stabilizer.
Then the roots of $(n\mapsto F_i(g_i^n))$ 
generate a finite index subgroup of the group
generated by the roots of $g_i$.
\end{lem}

\begin{pf}
Let $\alpha_{1,i},\ldots,\alpha_{d_i,i}$ be the roots of $g_i$; since
they are multiplicatively independent, the algebra generated by the
functions $(n\mapsto \alpha_{j,i}^n)$, for $j=1,\ldots,d_i$, is
isomorphic, as described above,
to the algebra $\C[X_1^{\pm 1},\ldots,X_{d_i}^{\pm 1}]$. Writing
$F_i(X_1,\ldots,X_{d_i})=\sum_{(j_1,\ldots,j_{d_i})}
a_{i,(j_1,\ldots,j_{d_i})}X_1^{j_1}\cdots X_{d_i}^{j_{d_i}}$,
the  roots of $(n\mapsto F_i(g_i^n))$ are the numbers
$\alpha_{i,1}^{j_1}\cdots \alpha_{d_i}^{j_{d_i}}$
for which the corresponding coefficient $a_{i,(j_1,\ldots,j_{d_i})}$
does not vanish.
The rank of the group they generate is then the rank of the lattice
generated in ${\mathbf Z}^{d_i}$ by the exponents
$(j_1,\ldots,j_{d_i})$ corresponding to non zero coefficients
$a_{i,(j_1,\ldots,j_{d_i})}$.
Suppose now by contradiction that group generated by the roots of
$(n\mapsto F_i(g_i^n))$ has rank $d<d_i$; then the above mentioned
lattice has also rank $d$, so is generated by $d$ vectors
$(l_{1,k},\ldots, l_{{d_i},k})\in {\mathbf Z}^d$, for $k=1,\ldots,d$.
Then we can write $F_i(X_1,\ldots,X_{d_i})$ as
\begin{equation*}
F_i(X_1,\ldots,X_{d_i})=G_i(X_1^{l_{1,1}}\cdots
X_{d_i}^{l_{d_i,1}},\ldots,X_1^{l_{1,d}}\cdots X_{d_i}^{l_{d_i,d}})
\end{equation*}
for some Laurent polynomial
$G_i(T_1,\ldots,T_d)\in{\mathbf C}[T_1^{\pm 1}\ldots, T_d^{\pm 1}]$
in $d<d_i$ variables.
Let  $(b_1,\ldots,b_{d_i})\in\Z^{d_i}\setminus\{0\}$ be a nonzero vector
which is orthogonal to all the vectors $(l_{1,k},\ldots, l_{d_i,k})$,
for $k=1,\ldots,d$; consider the algebraic subgroup of $\G_m^{d_i}$
formed by the elements of the form $(t^{b_1},\ldots,t^{b_{d_i}})$, with
$t\in\G_m$.
Translations with respect to this subgroup, $i.e.$ maps of
the form
$(X_1,\ldots,X_{d_i})\mapsto(t^{b_1}X_1,\ldots,t^{b_{d_i}}X_{d_i})$,
leave invariant the zero set
of $F_i$, which gives the desired contradiction.
\end{pf}

\begin{lem}
\label{findex}
In the previous notation, suppose that $\f_1$ divides $\f_2$ in the
ring of power series with roots in $\mathcal{U}$. Then the roots of $\f_2$
generate a finite index subgroup of $\mathcal{U}$.
\end{lem}

\begin{pf}
Using as usual the isomorphism between the ring of power series with
roots in $\mathcal{U}$ and the ring $\C[X_1^{\pm 1},\ldots,X_r^{\pm 1}]$,
the power sums $\f_i$ will correspond to Laurent polynomials
$G_i\in\C[X_1^{\pm 1},\ldots,X_r^{\pm 1}]$. By assumption, $G_1$
divides $G_2$ in the ring of Laurent polynomials. If, by contradiction,
the lemma were false, then  after applying a change of variables, we
could write $G_2$ as a Laurent polynomial in $(X_1,\ldots,X_{r-1})$,
while the variable $X_r$ would appear in $G_1$. Also, since the zero
set of $G_1$ has finite stabilizer, $G_1$ could not be of the form
$G_1(X_1,\ldots,X_r)=X_r^a \tilde{G_1}(X_1,\ldots,X_{r-1})$, for any
polynomial  $\tilde{G_1}(X_1,\ldots,X_{r-1}) $ independent of
$X_r$. But then  a divisibility relation of the form
$G_2(X_1,\ldots,X_{r-1})=G(X_1,\ldots,X_r)\cdot H(X_1,\ldots,X_r)$
could not be valid, and this contradiction finishes the proof.
\end{pf}

\begin{lem}
\label{D2finite}
Assume that each component of $D_2$ has finite stabilizer,
and  the hypotheses of the previous Lemma \ref{findex}.
Then the roots of $\f_1$ generate
a finite index subgroup of $\mathcal{U}$.
\end{lem}

\begin{pf}
In the  notation of the proof of Lemma \ref{findex},
we must prove that it is impossible that $G_1$
can be written such that one of the variables $X_1,\ldots,X_r$ is
omitted. This is due to the fact that otherwise $G_2$ will be
the product of $G_1$ and a Laurent polynomial
containing the omitted variable, so $F_2$ too would  be reducible,
and moreover one of its irreducible factors would omit an indeterminate,
up to automorphism of $\G_m^{d_2}$.
This is in contradiction with the fact that each component of
$D_2$ has trivial stabilizer,
so we have proved so far that the roots of $\f_1$ generate
a finite index subgroup of $\mathcal{U}$.
\end{pf}

Finally, we have obtained so far that both the roots of $g_1$
and those of $g_2$ generate
finite index subgroups of $\mathcal{U}$. So in particular $d_1=d_2=r$, and
the two tori $\G_1,\G_2$ are isomorphic, both having dimension $r$.

\begin{lem}
\label{D2trivial}
Assume the hypotheses of the two preceding Lemmas and that $D_2$ has
 trivial stabilizer.
Then the roots of $\f_1$ generate the whole group $\mathcal{U}$.
\end{lem}

\begin{pf}
By the previous lemma, the roots of $\f_1$ generate a finite index
subgroup of $\mathcal{U}$. Suppose by contradiction that such an index is
larger than one.
This means that, after changing generators of $\mathcal{U}$,
$\f_1$ can be written as a function of
$u_1^p,u_2,\ldots,u_r$, for some prime $p$, while in $\f_2$ cannot. In
terms of $F_1,F_2$, this implies that $F_2$ can be written as a
polynomial in $X_1^p,X_2,\ldots,X_{d_2}$, so its stabilizer would be
non-trivial.
\end{pf}

{\it Proof of Theorem \ref{with_multiplicity}}.
We now finish the proof. Recall that
$\f_i(n)$ can be written as $F_i(g_i^n)$, where $F(X_1,\ldots,X_r)=0$
is a polynomial equation for $D_i$ and also as
$G_i(u_1^n,\ldots,u_r^n)$, where $G_i$ is a Laurent polynomial and
$u_1,\ldots,u_r$ are multiplicatively independent.
Geometrically this means that we have isogenies
$\phi_i:\G_m^r\to\G_m^r$ defined by
\begin{equation*}
 g_i=\phi_i(u_1,\ldots,u_r),
\end{equation*}
such that $G_i=F_i(\phi_i)$. The fact that $\f_1$ divides $\f_2$ in the
ring of power sums  means that $G_1$ divides $G_2$, in the ring of
Laurent polynomials.
This means that
$\phi_1^* \I (D_1)\supset \phi_2^* \I(D_2)$.
Also, by the previous Lemma \ref{D2trivial},
$\deg\phi_1=1$, so we finish the proof by putting
$\phi=\phi_2\circ\phi_1^{-1}$.
\rightline{\it Q.E.D.}
\medskip

The following example shows that one cannot expect in general that the
morphism $\phi$ sends $g_1$ to $g_2$. 
\begin{ex}
\label{es}
{\rm
We take $k={\mathbf Q}, \OS=\Z$, $\G_1=\G_2=\G_m$, $D_1=D_2=1$ and
$g_1=2, g_2=-2$. We obtain that for even values of the exponent $n$,
$g_1^n=g_2^n$, so in particular the ideals $(g_1^n)^* \I (D_1)$ and
$(g_2^n)^* \I (D_2)$ coincide (both ideals are generated by the integer
$2^n-1$). Nevertheless there exists no morphism $\phi:\G_1\to \G_2$
sending $2\mapsto -2$ and satisfying
$\phi^* \I (D_2)\subset \I (D_1)$.
}
\end{ex}

Here is an example in which the divisor $D_2$ is reducible,
has trivial stabilizer, but its components have non trivial stabilizers,
so that hypothesis of Theorem \ref{with_multiplicity} is not satisfied.

\begin{ex}
\label{es2}
{\rm We take $k={\mathbf Q}, \OS=\Z$, $\G_1=\G_m$, $D_1=\{1\}$;
now put $\G_2=\G_m^2, D_2=\{1\}\times\G_m+\G_m\times \{1\}$,
so that $D_2=F_2^{-1}(0)$, for the polynomial $F_2(X_1,X_2)=(X_1-1)(X_2-1)$.
Choose $g_1=2, g_2=(2,3)$.
Clearly, condition \eqref{1.5} is satisfied for every $n$,
since it amounts to the fact that $2^n-1$ divides $(2^n-1)(3^n-1)$.
There exists no dominant map $\G_1\to\G_2$,
so the conclusion of Theorem \ref{with_multiplicity} fails. 
}
\end{ex}

Also, it may happens that $D_2$ has a non trivial stabilizer, while the
stabilizer of each of its component is trivial. In this case too, the
conclusion of Theorem \ref{with_multiplicity} can fail, as shown by the
following

\begin{ex}\label{es3}{\rm
Again, let us take $k={\mathbf Q}, \OS=\Z$, $\G_1=\G_m$, $D_1=\{1\}$ (so
$F_1(X)=X-1$). Then take $\G_2=\G_m, D_2=\{1, -1\}$, so $F_2(X)=X^2-1$
and  $\st (D_2)=\{\pm 1\}$.
Take $g_1=4$, $g_2=2$. For every positive integer $n$, $F_1(n)=F_2(n)$;
nevertheless, there is no integer $h>0$ and morphism $\phi$ such that
$\phi(g_1^h)=g_2^h$.
}
\end{ex}

Here is a general statement (the analogue of Theorem \ref{general}),
which essentially follows from the proof of Theorem
\ref{with_multiplicity}:

\begin{thm}\label{gene}
 Let $\OS$ be a ring of $S$-integers in a number field $k$.
Let $\G_1$ and $\G_2$ be linear tori, and let
$g_i\in\G_i(\OS)$  be elements generating Zariski-dense subgroups in
 $\G_i$ $(i=1,2)$.
Let $D_i$ be reduced divisors defined over $\OS$, with defining ideals
$\I(D_i)$. Suppose that for
infinitely many natural numbers $n\in\N$, the inclusion of ideals
\begin{equation*}
 (g_1^n)^* \I(D_1)\supset (g_2^n)^* \I(D_2)
\end{equation*}
holds. Then there exists a torus $\G_0$, dominant morphisms
 $\phi:\G_1\to \G_0$  and  $\psi:\G_2\to \G_0$ and a divisor $E$ on
 $\G_0$ such that $\phi(g_1)=\psi(g_2)$ and 
\begin{equation*}
 D_1\subset \phi^* E , \qquad D_2\supset \psi^* E.
\end{equation*}
\end{thm}

Let us show that the Examples \ref{es2}, \ref{es3}  can be treated by
Theorem \ref{gene}. In the situation of Example \ref{es2}, take $\G_0=\G_m$,
$\psi_1(X)=X$, $\psi_2(X,Y)=X$; in Example \ref{es3}, take $\G_0=\G_m,
\psi(X)=X^2, \phi(X)=X$.
\medskip

To prove Theorem \ref{gene} we need the following lemma, for which we
need a definition: we say that a power sum is {\it reduced} if in the
decomposition \eqref{powersum} one of the roots $\alpha_i$ is equal to
$1$. Notice that whenever $D$ is a divisor in a torus $\G=\G_m^r$, and
$g\in \G$ is a point generating a Zariski-dense subgroup, an equation
for $D$ can always be found of the form $F=0$, for a Laurent polynomial
$F\in k[X_1^{\pm 1},\ldots,X_r^{\pm 1}]$  such that the power sum
$n\mapsto \f(n):=F(g^n)$ is reduced. Actually, we can also take $F$ to
be a polynomial in $k[X_1,\ldots,X_r]$.

\begin{lem}\label{x}
Let $\f_1,\f_2,{\bf g}$ be power sum, such that the group generated by
all their roots is torsion free and $\f_1$ is reduced. Suppose that
$\f_2$ factors as $\f_2=\f_1\cdot {\bf g}$. Then every root of $\f_1$
has finite index in the group generated by the roots of $\f_2$.
\end{lem}

\begin{pf}
Suppose by contradiction  that one root of $\f_1$,  say $\gamma$, does
 not have finite index in the group generated by the roots of
$\f_2$. Then the group generated by the roots of $\f_1,\f_2$
and ${\bf g}$ admits a basis of the form $\gamma_1,\ldots,\gamma_r$,
where $\gamma_r^d=\gamma$, for a suitable $d>1$,
and $\gamma_1,\ldots,\gamma_{r-1}$ generate a group containing
the roots of $\f_2$. After the usual identification between
power sums with roots in a torsion-free rank $r$ group
and Laurent polynomials in $r$ indeterminates,
we can write $\f_2(n)=F_2(\gamma_1^n,\ldots,\gamma_{r-1}^n)$,
$\f_1(n)=F_1(\gamma_1^n,\ldots,\gamma_r^n)$ and
${\bf g}(n)=G(\gamma_1^n,\ldots,\gamma_r^n)$,
for Laurent polynomials $F_1,F_2,G\in k[X_1^{\pm 1},\ldots,X_r^{\pm 1}]$,
where $F_2$ does not depend on $X_r$,
while $F_1$ does. From the factorization $\f_2=\f_1\cdot{\bf g}$
follows the corresponding factorization $F_2=F_1\cdot G$.
Since $\f_1$ is assumed to be reduced, $F_1$ is not of
the form $X_2^k\cdot \tilde{F_1}$, for any $k\neq \Z$
and $\tilde{F_1}\in k[X_1^{\pm 1},\ldots, X_{r-1}]$.
Hence $F_1$ is not invertible in $k(X_1,\ldots,X_{r-1})[X_r^{\pm 1}]$,
so it is impossible that $F_2$ omits
the indeterminate $X_r$, finishing the proof.
\end{pf}

{\it Proof of Theorem \ref{gene}}. Let, as usual, $\G_i=\G_m^{d_i}$,
$g_i=(\alpha_{1,i},\ldots,\alpha_{d_i,i})$. Let $F_i=0$ be equations for
$D_i$, such that the power sum $n\mapsto F_1(g_1^n)=:\f_1(n)$ is
reduced.  The inclusion of ideals in our assumption means that $\f_1(n)$
divides $\f_2(n)$ for infinitely many integers $n$. As in the previous
proofs, we can reduce to the case when the roots of the two power sums
together generate a torsino-free multiplicative group (this is obtained
by considering separately the  arithmetic progressions $n\mapsto qn+r$,
for a suitable $q>1$, $r\in\{0,\ldots, q-1\}$, and applying the
assumptions to each such arithmetic progression).  Hence we can apply
Lemma \ref{cz}, which provides a power sum $n\mapsto {\bf g}(n)$ such
that identically $\f_2=\f_1\cdot{\bf g}$. By Lemma \ref{x}, the roots of
$\f_1$ have finite index in the group generated by the roots of
$\f_2$. Let $q$ be the minimal common multiple of such indices. Let $r$
be the rank of the group $\Gamma$ generated by the roots of $\f_1$ and
put $\G_0=\G_m^r$. Then the group $\Gamma^d:=\{\gamma^d\, :\,
\gamma\in\Gamma\}$ embeds both in the group generated by
$\alpha_{1,1},\ldots,\alpha_{d_1,1}$ and in the group generated by
$\alpha_{1,2},\ldots,\alpha_{d_2,2}$. Let $\gamma_1,\ldots,\gamma_r$ be
a basis of $\Gamma$.
The embeddings of $\Gamma$ in the two mentioned groups correspond  to
dominant morphisms $\phi: \G_1\to\G_0$ and $\psi:\G_2\to\G_0$ with
$\phi(g_1)=\psi(g_2)=(\gamma_1,\ldots,\gamma_r)$.  Let us write
$\f_1(qn)=F_0(\gamma_1,\ldots,\gamma_r)$ for a Laurent polynomial
$F_0\in k[X_1^{\pm 1},\ldots, X_r^{\pm 1}]$. Then, putting
$E=F_0^{-1}(0)$ we have $D_1\subset\phi^* E,\, D_2\supset\psi^*(E)$ as
wanted.
\rightline{\it Q.E.D.}

\section{Arithmetic support problem.}

In this section we study the unicity problem only with
supports in the assumption of Theorem \ref{with_multiplicity}.
We need some technical hypothesis on the divisors $D_i$.
We prove the following.

\begin{prop}\label{no_multiplicity}
Let $\G_1,\G_2,g_1,g_2$ be as above, let $D_1,D_2$ be irreducible
divisors such that $D_1$ contains the origin of $\G_1$ and $D_2$ does
not contain any translate of a positive dimensional sub-torus. Suppose
that for every sufficiently large integer $n$ the inclusion
\begin{equation}
\label{radicals}
 \supp (g_1^n)^* \I (D_1)\subset \supp (g_2^n)^* \I (D_2) 
\end{equation}
holds. Then there exists a dominant morphism $\phi: \G_1\to\G_2$,
 defined over $k$, and an integer $h \geq 1$ such that 
$\phi(g_1)=g_2^h$.
\end{prop}

We note that the condition for $D_2$ to contain no translate of any
positive dimensional sub-torus is much stronger than
the (necessary) condition that its stabilizer be trivial.
We do not know whether the latter suffices. Also,
we do not obtain any relation between $\phi^* D_2$ and $D_1$.
In the case we have equality of supports, we obtain the stronger
conclusion that $\phi$ is \'etale.

The above results will be reduced to a theorem of Barsky, B\'ezivin and
Schinzel \cite{bbs}.
We state as a lemma a particular case of Theorem 1 of \cite{bbs}:

\begin{lem}\label{Lemmabbs}
Let $k$ be a number field, let $\OS\subset k$ be a ring of $S$-integers,
let $\alpha_1,\ldots,\alpha_{d_1}\in \OS^*$ be multiplicatively independent
units in $\OS$,
and let $\beta_1,\ldots,\beta_{d_2}\in \OS^*$ be units of $\OS$.
Let $F_1(X_1,\ldots,X_{d_1})$ and $F_2(Y_1,\ldots,Y_{d_2})$ be
polynomials.
Assume that $F_1(1,\ldots,1)=0$ and that the
equation $F_2(y_1,\ldots,y_{d_2})=0$ has only finitely many solutions
in the roots of unity $y_1,\ldots,y_{d_2}$.
If for all large integers $n$ 
\begin{equation}
\label{supprel}
\supp\, \I (F_1(\alpha_1^n,\ldots,\alpha_{d_1}^n)) \subset
\supp\, \I (F_2(\beta_1^n,\ldots,\beta_{d_2}^n)),
\end{equation}
then there exists a positive integer $h$ such that
 $\beta_1^h,\ldots,\beta_{d_2}^h$ belongs to the multiplicative group
 generated by $\alpha_1,\ldots,\alpha_{d_1}$.
\end{lem}

{\it Proof of Proposition \ref{no_multiplicity}}.
Write
$g_1=(\alpha_1,\ldots,\alpha_{d_1})\in\G_m^{d_1}$ and
$g_2=(\beta_1,\ldots,\beta_{d_2})$. The divisors $D_i\subset\G_m^{d_i}$
will be defined in $\G_m^{d_i}$ by the equation
$F_i(X_1,\ldots,X_{d_i})=0$, where $F_i$ are irreducible polynomials.
Suppose the hypotheses of Theorem \ref{no_multiplicity} are
satisfied. In particular, $\alpha_1,\ldots,\alpha_{d_1}$ are
multiplicatively independent, $\beta_1,\ldots, \beta_{d_2}$ are also
multiplicatively independent and $F_1(1,\ldots,1)=0 $. The fact that the
divisor $D_2$ contains no translate of positive dimensional sub-tori
implies, by a theorem of Laurent (\cite{lau}, previously a conjecture of
Lang), that it contains only finitely many points whose coordinates are
roots of unity. So the equation $F_2(y_1,\ldots,y_{d_2})=0$ has only
finitely many solutions in roots of unity, as required in Lemma
\ref{Lemmabbs}.
The hypothesis that
$\supp\, (g_1^n)^* \I (D_1) \supset \supp\, (g_2^n)^* \I (D_2)$
is equivalent to condition \eqref{supprel},
so Lemma \ref{Lemmabbs} applies
and gives the existence of an integer
$h \geq 1$ and a $d_2\times d_1$ matrix
$(a_{ij})_{1\le i\le d_2; 1\le j\le d_1}$ such that
\begin{equation*}
 \beta_i^h=\prod_{j=1}^{d_1}\ \alpha_j^{a_{ij}}. 
\end{equation*}
Then, defining $\phi:\G_m^{d_1}\to\G_m^{d_2}$ by sending
$(x_1,\ldots,x_{d_1})\mapsto (\prod_{j}x_j^{a_{1j}},\ldots,\prod_{j}
x_j^{a_{d_2 j}}) $, we obtain the desired conclusion.
\hfill{\it Q.E.D.}

\bigskip
\baselineskip=12pt
\rightline{Dipartimento di Matematica e Informatica}
\rightline{University of Udine}
\rightline{Via delle Scienze, 206 - 33100 Udine}
\rightline{e-mail: pietro.corvaja@dimi.uniud.it}
\bigskip

\baselineskip=12pt
\rightline{Graduate School of Mathematical Sciences}
\rightline{The University of Tokyo}
\rightline{Komaba, Meguro,Tokyo 153-8914}
\rightline{e-mail: noguchi@ms.u-tokyo.ac.jp}

\medskip

\begin{thebibliography}{99}
\setlength{\itemsep}{-3pt}
\bibitem{bbs} Barsky D., B\'ezivin, J-P. and Schinzel, A., Une
 caract\'erisation arithm\'etique de suites r\'ecurrentes lin\'eaires,
 J. Reine Angew. Math. {\bf 494} (1998), 73-84.
\bibitem{crs}
Corrales-Rodorig\'a\~nez, C. and Schoof, R.,
The support problem and its elliptic analogue,
J. Number Theory {\bf 64} (1997), 276-290.
\bibitem{czIndag} Corvaja, P. and Zannier, U., Diophantine equations with
 power sums and universal Hilbert sets, Indag. Math. {\bf N.S. 9(3)}
 (1998), 317-332.
\bibitem{cz}Corvaja, P. and Zannier, U., Finiteness of integral values for
 the ratio of two linear recurrences, Invent. Math. {\bf 149} (2002),
 431-451.
\bibitem{lau} Laurent, M., Equations diophantiennes exponentielles,
Invent. Math. {\bf 78} (1984), 299-327.
\bibitem{n81}
Noguchi, J.,
Lemma on logarithmic derivatives and holomorphic curves in algebraic
 varieties, Nagoya Math.\ J.\ {\bf 83} (1981), 213-233.
\bibitem{n98}
Noguchi, J.,
On holomorphic curves in semi-Abelian varieties,
Math.\ Z. {\bf 228} (1998), 713-721.
\bibitem{no84}
Noguchi, J. and Ochiai, T.,
Geometric Function Theory in Several Complex Variables,
Japanese edition, Iwanami, Tokyo, 1984;
English Translation, Transl.\ Math.\ Mono.\ {\bf 80},
Amer.\ Math.\ Soc., Providence, Rhode Island,
1990.
\bibitem{nwy02}
Noguchi, J., Winkelmann, J. and Yamanoi, K.,
The second main theorem for holomorphic curves into semi-Abelian
varieties,
Acta Math.\ {\bf 188} no.1 (2002), 129-161.
\bibitem{nwy08}
Noguchi, J., Winkelmann, J. and Yamanoi, K.,
The second main theorem for holomorphic curves into semi-Abelian
varieties II,  Forum Math.{\bf 20} (2008), 469-503.
\bibitem{vdP} {van der Poorten, A. J.},
{Solution de la conjecture de {P}isot sur le quotient de {H}adamard de
deux fractions rationnelles},
{C.\ R.\ Acad.\ Sci.\  S\'er.\ I Math.} 
 {\bf 306}   (1988), 97--102.
\bibitem{rum} Rumely, R., Notes on van der Poorten's proof of the
Hadamard Quotient Theorem: Part I, II, S\'eminaire de Th\'eorie des
Nombres de Paris 1986-87, pp.\ 349--382, pp.\ 383--409, 
Progr.\ Math.\ {\bf 75}, Birkh\"auser, Boston-Basel, 1988. 
\bibitem{S}  Schmidt, W. M.,
Linear recurrence sequences, 
in Diophantine Approximations, pp. 171-247, Lecture Notes in Math.\ {\bf 1819},
Springer, Berlin, 2003.
\bibitem{y04}
Yamanoi, K.,
Holomorphic curves in abelian varieties and intersection
with higher codimensional subvarieties, 
Forum Math.\ {\bf 16} (2004), 749-788.
\end{thebibliography}
\end{document}